\documentclass[12pt,reqno]{amsart}

\usepackage[dvips]{graphicx} 

\usepackage{
hyperref
}

\usepackage{breakurl}   

\usepackage{epigraph}

\theoremstyle{definition}

\numberwithin{equation}{section}

\title{Is pluralism in the history of mathematics possible?}

\author[J. B.]{Jacques Bair}\address{HEC-ULG, University of Liege,
  4000 Belgium}\email{j.bair@ulg.ac.be}

\author[A. B.]{Alexandre Borovik} \address{School of Mathematics,
  University of Manchester, Oxford Street, Manchester, M13 9PL, United
  Kingdom} \email{alexandre@borovik.net}

\author[V. K.]{Vladimir Kanovei} \address{IITP RAS, Moscow,
  Russia}\email{kanovei@googlemail.com}

\author[M. K.]{Mikhail G. Katz}\address{Department of Mathematics, Bar
  Ilan University, Ramat Gan 5290002
  Israel}\email{katzmik@math.biu.ac.il}

\author[S. K.]{Semen S. Kutateladze}\address{Sobolev Institute of
  Mathematics, Novosibirsk State University, Russia}
\email{sskut@math.nsc.ru}

\author[S. S.]{Sam Sanders} \address{Department of Philosophy 2, RUB
  Bochum, Bochum, Germany
  \url{http://sasander.wix.com/academic}}\email{sasander@me.com}

\author[D. S.]{David Sherry} \address{Department of Philosophy,
  Northern Arizona University, Flagstaff, AZ 86011, US}
\email{David.Sherry@nau.edu}

\author[M. U.]{Monica Ugaglia} \address{Il Gallo Silvestre, Localit\`a
  Collina 38, Montecassiano, Italy}\email{monica.ugaglia@gmail.com}

\author[M. V.]{Mark van Atten} \address{Archives Husserl (CNRS / ENS),
  45 rue d’Ulm, 75005 Paris, France} \email{mvanatten@ens.fr}

\begin{document}

\thispagestyle{empty}


\begin{abstract}
Leibniz scholarship is currently an area of lively debate.  We respond
to some recent criticisms by Archibald et al.
\end{abstract}

\maketitle



This is a brief response to the article ``A Question of Fundamental
Methodology: Reply to Mikhail Katz and his coauthors,'' by Archibald
et al. in \emph{The Mathematical Intelligencer} \cite{Ar22b}.  The
article by Archibald et al.~was in reaction both to our earlier
article ``Two-track depictions of Leibniz's fictions'' \cite{22b} in
the same journal, and to other work of ours.

We have argued that, in addition to procedures that can be adequately
described in purely Archimedean settings, Leibniz (as well as Cauchy
and others) used procedures that exploit genuine infinitesimals, that
is, what to him were mathematical entities.

Archibald's coauthors Arthur and Rabouin have argued that the term
``infinitesimal'' as used by Leibniz does not refer to a mathematical
entity, and is, rather, stenography for exhaustion-type arguments in
the style of Archimedes.  We have compared the two approaches in
\cite{22b} and, in particular, presented evidence for our
interpretation.

Archibald et al.~make a number of false claims concerning both
\cite{22b} and other publications of ours.  Lack of space prevents us
from responding in full.  A more detailed response appears at
\cite{22a}.

In closing, it is ironic that Archibald et al.~should claim that
\begin{enumerate}\item[]
[O]ver the years, it became clearer and clearer that our interlocutors
do not care much about rational discussion and scientific dialogue
from \emph{different perspectives}, {\ldots} The latest example of
that approach is provided by a paper {\ldots}~``Two-Track Depictions
of Leibniz’s Fictions.'' \cite[p.\;2]{Ar22b} (emphasis on ``different
perspectives'' added)
\end{enumerate}
``Two-track depictions'' \cite{22b} is devoted specifically to making
explicit a pair of \emph{different perspectives} on Leibniz's
calculus, so as to stimulate rational discussion and scientific
dialogue.

Archibald et al.~do little to clarify the ``Question of Fundamental
Methodology'' of their title, namely that the history of mathematics,
like mathematics itself, could benefit from a plurality of approaches.


\begin{thebibliography}{AII}


\bibitem{Ar22b} Archibald, Tom; Arthur, Richard T. W.; Ferraro,
  Giovanni; Gray, Jeremy; Jesseph, Douglas; L\"utzen, Jesper; Panza,
  Marco; Rabouin, David; Schubring, Gert.\, A Question of Fundamental
  Methodology: Reply to Mikhail Katz and His Coauthors.  \emph{The
    Mathematical Intelligencer} \textbf{44} (2022), no.\;4, 360--363.
  \url{https://doi.org/10.1007/s00283-022-10217-7} \MR{4526076}

\bibitem{22a} Bair, J.; Borovik, A.; Kanovei, V.; Katz, M.;
  Kutateladze, S.; Sanders, S.; Sherry, D.; Ugaglia, M.\, Historical
  infinitesimalists and modern historiography of infinitesimals.
  \emph{Antiquitates Mathematicae} \textbf{16} (2022), 189--257.
  \url{https://arxiv.org/abs/2210.14504},
  \url{https://doi.org/10.14708/am.v16i1.7169}  

\bibitem{22b} Katz, M.; Kuhlemann, K.; Sherry, D.; Ugaglia, M.; van
  Atten, M.\, Two-track depictions of Leibniz's fictions.  \emph{The
    Mathematical Intelligencer} \textbf{44} (2022), no.\;3, 261--266.
  \url{https://doi.org/10.1007/s00283-021-10140-3},
  \url{https://arxiv.org/abs/2111.00922} \MR{4480193}


\end{thebibliography}
\end{document}